\documentclass[11pt,twoside,a4paper]{article}
\usepackage{myartmods}
\usepackage[utf8]{inputenc}
\usepackage{datetime}\usdate
\usepackage[english]{babel}

\usepackage{amsmath}
\usepackage{amsfonts} 
\usepackage{amssymb}
\usepackage{mathtools}
\usepackage{multirow,bigdelim}

\usepackage{graphicx}
\usepackage[caption=false]{subfig}
\usepackage{tabularx}
\usepackage{booktabs} 
\usepackage{authblk}
\usepackage{color,soul} 
\soulregister\cite7
\soulregister\ref7
\soulregister\eqref7
\usepackage[usenames,dvipsnames]{xcolor}

\usepackage[colorlinks]{hyperref}

\numberwithin{equation}{section}
\numberwithin{figure}{section}
\numberwithin{table}{section}

\catcode`@=11
\def\HI{\text{HI}}
\def\tHI{\text{tHI}}
\def\minim{\mathop{\operator@font minimize}}
\def\minimize#1{{\displaystyle\minim_{#1}}}
\def\subject{\mathop{\operator@font subject\ to}}
\catcode`@=12

\begin{document}
\title      {Increased accuracy of planning tools for optimization of dynamic multileaf collimator delivery of radiotherapy through reformulated objective functions}
\author[1,2]{Lovisa Engberg\thanks{Corresponding author: loven140@kth.se}}
\author[2]  {Kjell Eriksson}
\author[1]  {Anders Forsgren}
\affil[1]   {Optimization and Systems Theory, Department of Mathematics, \protect\\KTH Royal Institute of 
             Technology, Stockholm SE-100 44, Sweden}
\affil[2]   {RaySearch Laboratories, Sveav\"{a}gen 44, Stockholm SE-103 65, Sweden}
\date       {Revised manuscript\\April 19, 2018}
\markboth   {Increased accuracy of planning tools through reformulated objective functions}{Increased accuracy of planning tools through reformulated objective functions}
\maketitle\thispagestyle{empty}

\begin{abstract} 
The purpose of this study is to examine in a clinical setting a novel formulation of objective functions for intensity-modulated radiotherapy treatment plan multicriteria optimization (MCO) that we suggested in a recent study. The proposed objective functions are extended with dynamic multileaf collimator (DMLC) delivery constraints from the literature, and a tailored interior point method is described to efficiently solve the resulting optimization formulation. In a numerical planning study involving three patient cases, DMLC plans Pareto optimal to the MCO formulation with the proposed objective functions are generated. Evaluated based on pre-defined plan quality indices, these DMLC plans are compared to conventionally generated DMLC plans. Comparable or superior plan quality is observed. Supported by these results, the proposed objective functions are argued to have a potential to streamline the planning process, since they are designed to overcome the methodological shortcomings associated with the conventional penalty-based objective functions assumed to cause the current need for time-consuming trial-and-error parameter tuning. In particular, the increased accuracy of the planning tools imposed by the proposed objective functions has the potential to make the planning process less complicated. These conclusions position the proposed formulation as an alternative to existing methods for automated planning. 

\end{abstract}

\section{Introduction}
Optimization of intensity-modulated radiation therapy (IMRT) treatment plans using the widely used penalty-based objective functions results in plans of high quality in a majority of cases (see, e.g., \cite{bokrantz2013,tol2015b} for references on the conventional penalty functions, and \cite{bortfeld1999} for a general discussion on IMRT treatment planning based on these functions). The path towards satisfactory clinical quality, however, is not always straight. Due to an implicit and not so clear relationship between the conventional objective functions and measures of plan quality, several re-optimizations with manually revised parameters are often needed---if not to improve the plan quality, then to explore local variations in order to feel convinced on the optimality.

Limitations of the current planning tools are also reflected in the existence of literature on methods for detection and/or post-processing of suboptimal plans. For example, Appenzoller et al.~\cite{appenzoller2012} suggest a method to predict achievable dose-volume histograms (DVHs) to be used as reference, and Fredriksson \cite{fredriksson2012} describes a technique to reduce the dose delivered to healthy tissue subsequent to optimization. Fredriksson successfully improves in that sense an already optimized plan while maintaining target coverage. 

The large amount of resources spent on treatment planning has motivated the development of strategies to automate the planning process. Automating strategies include machine learning, with the aim of predicting planning related parameters based on a database with previously treated patients; and computerization of planning schemes, often with the aim of mimicking the planning steps of the human planner by enclosing the optimization in an outer automated fine-tuning loop. For example, McIntosh and Purdie \cite{mcintosh2017} and Shiraishi and Moore \cite{shiraishi2016} present methods to predict an achievable dose distribution to possibly input in a subsequent automated plan reconstruction phase, and Wu et al.~\cite{wu2011} and Appenzoller et al.~\cite{appenzoller2012} suggest methods to predict parts of or entire achievable DVHs to support the planner's desicions on planning parameters. Some automating strategies have already been implemented into clinical practice through commercial treatment planning systems (TPSs). Yuan et al.~\cite{yuan2012} describe the machine-learning concepts of one commercially available DVH prediction method, which has been later evaluated, e.g., in the study by Tol et al.~\cite{tol2015}. Gintz et al.~\cite{gintz2016} describe and evaluate a clinically available computerization strategy to mimic a planning scheme that can handle a general patient case. 

The common approach of many automating strategies is to remove or largely reduce human interaction within the planning process. A risk associated with such an approach is ``automation surprises'', causing the user to feel out of control and leading to unnecessary cognitive load \cite{baxter2012}. In our research, we take a different perspective on automated treatment planning. We argue that the demand for automation is partly due to methodological shortcomings in the conventional penalty-based formulation, and that resolving these issues by reformulating the objective functions have a similar potential as other automating strategies to streamline the treatment planning process. A qualified formulation of objective functions has to offer a stronger connection to plan quality measures and lead to tractable (e.g., convex) optimization problems, yet without compromising performance---the plan quality achieved must be at least comparable to that achieved from conventional planning. A streamlining effect would then be expected as a result of eliminated or reduced trial and error due to increased accuracy of the planning tools provided in the TPS.

In a recent fluence-based study \cite{engberg2017}, we suggest a novel formulation of objective functions that possesses the above-mentioned properties to streamline the planning process. We abandon the implicit penalty-function based paradigm to obtain a stronger connection to commonly used plan quality measures, and use approximations in order to achieve a convex optimization problem. The results of our numerical planning study indicate an ability of the proposed formulation to perform equivalently or better (in terms of given plan quality measures) as compared to the conventional formulation in the given fluence-based setting. As to the study presented in this manuscript, the purpose is to strengthen our conclusions regarding the streamlining potential of the proposed formulation by evaluation in a more clinical context. To this end, deliverability constraints from the literature are added to handle the delivery with the dynamic multileaf collimator (DMLC) technique. This extension enables clinically accurate dose computations, and thus comparison with DMLC plans conventionally generated using clinical TPSs.

\section{Method}\label{sec:Methods}
The method of this study is an extension of the objective functions introduced in our previous fluence-based study \cite{engberg2017} with the linear DMLC deliverability constraints suggested by Papp and Unkelbach \cite{papp2014}; and the design of a tailored interior point method to efficiently solve the resulting MCO formulation. In Section~\ref{sec:Methods:ProposedFormulation}, we state and give a motivation for the proposed objective functions, and review the deliverbility constraints by Papp and Unkelbach in Section~\ref{sec:Methods:Deliverability}. Finally, in Section~\ref{sec:Methods:Solvability}, we discuss the matter of solving the entire MCO formulation. In particular, we demonstrate the linear-algebraic steps towards achieving the tailored interior point method, which are analogous to the derivation of the method described in \cite{engberg2017} for the deliverability-relaxed MCO formulation.

\subsection{Formulation of MCO objective functions}\label{sec:Methods:ProposedFormulation}
The plan quality indices that we consider are those which can be deduced from the dose-volume histograms (DVHs) of the dose distribution. Specifically, a \emph{dose-at-volume} index here refers to the dose level at a given volume fraction in the DVH of a region of interest (ROI), such as a planning target volume (PTV) or an organ at risk (OAR), thus corresponding to the minimum dose received by this fraction of the particular ROI. In treatment plan optimization, doses-at-volume for possibly several volume fractions of the ROIs are to be either minimized or maximized, and/or kept below or above given dose bounds, to attain highest possible (quantitative) plan quality and fulfil clinical requirements. This definition of the general planning goal is in line with the ``Plan Quality Metric'' designed by Nelms et al.~\cite{nelms2012} to unambiguously quantify and compare the plan quality of different plans. The conventional way of mathematically formulating the planning goal is by objective functions that penalize any deviation of the DVH curves from user-specified reference points (see, e.g., \cite{bokrantz2013,tol2015b}). Since the connection between such penalties and the actual dose-at-volume is implicit and not always clear, the optimization operation might behave unexpectedly, causing a need for re-optimization and tedious parameter tuning. On the contrary, the ideal way of formulating the planning goal is, naturally, by objective functions that explicitly quantify the doses-at-volume in question. A consequence of such a choice would be an intractable optimization problem which cannot be solved to proven optimality: the function quantifying dose-at-volume is both non-differentiable and non-convex. The idea behind our proposed formulation is to adopt the ideal explicit structure, by which the penalty-function paradigm is abandoned, while replacing the dose-at-volume objective functions by convex over- and concave underestimates. We use the mean-tail-dose functions introduced by Romeijn et al.~\cite{romeijn2006} to obtain such estimates. A brief review of the proposed MCO formulation follows here; the full formulation is stated in Appendix~\ref{app:propForm} and a detailed description is given in our previous paper \cite{engberg2017}. 

Let $d$ be a distribution of dose inside the patient geometry. Denote by $D = D(d;\,v,s)$ the function quantifying the dose-at-volume associated with the volume fraction $v$ of the ROI $s$, and by $D^+ = D^+(d;\,v,s)$ and $D^- = D^-(d;\,v,s)$ the functions quantifying the upper and lower mean-tail-dose, respectively. The proposed MCO formulation takes $D^+$ as objective function to approximate the minimization of $D$, and analogously the negative of $D^-$ to approximate the maximization. With $K$ doses-at-volume out of which $q$ should be minimized and $K-q$ maximized, the formulation is
\begin{equation}\label{eq:propForm}
\begin{aligned}
& \minimize{d} 
           && \mkern-10mu \mathmakebox[0pt][l]{\big[\,\xi_1, \cdots, \xi_q, -\xi_{q+1}, \cdots, -\xi_K\,\big]^T} \\
& \subject && D^+(d;v_k,s_k) \leq \xi_k \leq u_k, && \hat{l}_k \leq \xi_k, && k = 1,\ldots,q,       \\ 
&          && D^-(d;v_k,s_k) \geq \xi_k \geq l_k, && \hat{u}_k \geq \xi_k, && k = q\!+\!1,\ldots,K, \\
&          && \mathmakebox[0pt][l]{d \text{ deliverable dose distribution,}}
\end{aligned}
\end{equation} 
where auxiliary variables $\xi_k$ are introduced for clarity. Hard upper and lower dose bounds $u_k$ and $l_k$ are imposed on the mean-tail-doses via their auxiliary variable, and should be used to exclude clinically irrelevant dose distributions. Lower and upper dose bounds $\hat{l}_k$ and $\hat{u}_k$ are imposed only on the auxiliary variables, giving the effect of removed incentives to further minimize or maximize the mean-tail-doses; thus, $\hat{l}_k$ and $\hat{u}_k$ should be used to specify clinically satisfactory (or utopian) dose levels. Alternatively, excluding the auxiliary variable from the multicriteria objective removes all optimizing incentive and leaves only the hard dose bound for that mean-tail-dose. Notice that the ideal formulation is obtained by simply replacing $D^+$ and $D^-$ with $D$. Since the upper and lower mean-tail-doses over- and underestimate the dose-at-volume, i.e.,
\[
	D^+(d;v,s) \geq D(d;v,s) \geq D^-(d;v,s), \quad \text{for all } d, 
\]
the dose-at-volume is properly controlled and guaranteed to satisfy the hard dose bounds imposed on the mean-tail-doses. Regarding conservativeness of using the clinical dose bounds aimed for $D$ as upper and lower dose bounds for $D^+$ and $D^-$, this is briefly discussed in Section~\ref{sec:Discussion}. 

By convexity of $D^+$ and concavity of $D^-$, the formulation in \eqref{eq:propForm} is a convex MCO formulation provided that the intended set of deliverable dose distributions is convex (with respect to the delivery parameters). Without compromising convexity, \eqref{eq:propForm} can be extended with objective functions to minimize the maximum dose, maximize the minimum dose, or in any direction optimize the average dose. These extensions are included in the full formulation given in Appendix~\ref{app:propForm}.

\subsection{Modelling of DMLC deliverability}\label{sec:Methods:Deliverability}
We use the formulation by Papp and Unkelbach \cite{papp2014} to mathematically model the deliverability requirement in \eqref{eq:propForm} for ``sliding window'' DMLC delivery. DMLC delivery allows the collimator leaves to move during irradiation. By restricting to unidirectional sweeps (sliding windows), Papp and Unkelbach were able to express the resultant dose distribution using only linear transformations. This was achieved by representing the sweeps as the sequences of time of arrival and departure of the leaves at fixed positions (bixels) along the sweeping direction, rather than using control-point based representations where the trajectories are given by the sequence of the positions of the leaves at fixed points in time. Control-point based representations require in general a non-convex mapping to transform into dose distribution (see, e.g., \cite{unkelbach2015}). A brief review of the formulation by Papp and Unkelbach, with some supplementary assumptions, follows here. 

Denote by $B$ the number of beam angles, by $N$ the number of leaf pairs of the MLC, and by $J$ the number of bixels swept by each leaf pair, so that the fluence map for each beam is defined on a rectangular $N \times J$ grid of bixels. The trajectory of the $n$th leaf pair in the $b$th beam is defined as the sequences of departure times $l_{b,n,j}$ and $r_{b,n,j}$, for all bixels $j$, when the trailing (WLOG, the left) and the leading (right) leaf begin traversing each bixel. Assuming that all leaves require a constant time $\Delta t$ to traverse a bixel, the conditions for a feasible trajectory become
\begin{subequations}\label{eq:condOfExistence}
\begin{alignat}{2}
  r_{b,n,j} + \Delta t & \leq r_{b,n,j+1},   && \quad j = 1,\ldots,J\!-\!1, \label{eq:condOfExistence-a} \\
  l_{b,n,j} + \Delta t & \leq l_{b,n,j+1},   && \quad j = 1,\ldots,J\!-\!1, \label{eq:condOfExistence-b} \\
  r_{b,n,j}            & \leq l_{b,n,j},     && \quad j = 1,\ldots,J,       \label{eq:condOfExistence-c}
\end{alignat}
\end{subequations}
where \eqref{eq:condOfExistence-a} and \eqref{eq:condOfExistence-b} prevent a leaf to depart from bixel $j+1$ before arriving there, and \eqref{eq:condOfExistence-c} prevents a trailing leaf to pass a leading leaf. The effect of assuming constant traversing time is a halved number of DMLC variables as compared to the original formulation by Papp and Unkelbach, who also keep track of arrival times. The assumption reduces the set of feasible trajectories but, most importantly, has no impact on the set of deliverable dose distributions. Further restrictions to the trajectories can be made depending on limitations of the delivering treatment machine or other resources. For example, Papp and Unkelbach outline how to handle interdigitation constraints. The machine used in this study allows interdigitation, but requires the leading and trailing leaves to always be separated by a minimum gap corresponding to a fraction $\rho \in (0,1)$ of the bixel width. The leading leaf must thus depart $\rho\Delta t$ time units before the trailing leaf arrives, 
\begin{subequations}\label{eq:minGap}
\begin{alignat}{2}
  r_{b,n,j+1} & \leq l_{b,n,j} + \Delta t - \rho\Delta t,                \nonumber\\
              & =    l_{b,n,j} + (1-\rho)\Delta t,       && \quad j = 1,\ldots,J\!-\!1, \label{eq:minGap-a}\\
  r_{b,n,1}   & \leq l_{b,n,1} - \rho\Delta t.                           \label{eq:minGap-b}
\end{alignat}
\end{subequations}
Here, \eqref{eq:minGap-b} handles the case with the left-most bixel. Notice that the minimum gap constraints \eqref{eq:minGap} make \eqref{eq:condOfExistence-c} redundant. We also require the entire delivery to be completed within some maximum treatment time $T_{\max}$. Since the delivery of beam $b$ is completed when all trailing leaves have traversed the right-most bixel $J$, this requirement is formulated by the linear constraints
\begin{equation}\label{eq:maxTime}
\begin{aligned}
  l_{b,n,J} + \Delta t & \leq T_b,      && \quad n = 1,\ldots,N \text{ and } b = 1,\ldots,B, \\
  \sum_{\forall b} T_b & \leq T_{\max},                                                      
\end{aligned}
\end{equation}
where $T_b$ will be regarded as the beam-on time for beam $b$. 

The dose distribution is determined by the time of exposure to radiation of the bixels. As observed by Papp and Unkelbach, the exposure can be expressed linearly: the full exposure of the bixels traversed by leaf pair $n$ in beam $b$ equals the time intervals between the departure of the trailing and leading leaves, 
\[\begin{aligned}
	l_{b,n,j} - r_{b,n,j}, && j = 1,\ldots,J.
\end{aligned}\]
We slightly extend this model by accounting for partial exposure due to MLC leakage during the remaining portion of the beam-on time $T_b$, 
\[\begin{aligned}
	\tau \left( T_b - \left(l_{b,n,j} - r_{b,n,j}\right) \right), && j = 1,\ldots,J,
\end{aligned}\]
where $\tau \in (0,1)$ denotes the transmission factor. Assuming constant dose rate $\delta$, the beamlet weights of the fluence map are then given by
\[\begin{aligned}
    \delta\big(l_{b,n,j}-r_{b,n,j} + \tau\big(T_b - (l_{b,n,j}-r_{b,n,j})\big)\big), && j = 1,\ldots,J, 
\end{aligned}\]
and the dose distribution $d$ (now discretized into voxels) is obtained from multiplication with the dose deposition matrix $P$,
\begin{equation}\label{eq:doseComp}
\begin{aligned}
    d = \delta P \big(l-r + \tau\big(T - (l-r)\big)\big),
\end{aligned}
\end{equation}
where $l$ and $r$ are vectorizations of $l_{b,n,j}$ and $r_{b,n,j}$, and the vector $T$ is the appropriate vectorization of the beam-on times (i.e., the vector of $NJ$ repetitions of each $T_b$).

\subsection{On solving the proposed MCO formulation}\label{sec:Methods:Solvability}
A Pareto optimal solution to \eqref{eq:propForm} is obtained as the optimal solution to a weighted-sum instance, where the multiple objective functions are scaled by some weighting factors and accumulated into one. Different weighting factors lead to different Pareto optimal solutions. Rockafellar and Uryasev \cite{rockafellar1997} have shown that the conditional value-at-risk measure---the financial counterpart to mean-tail-dose---can be expressed as a linear program by the introduction of artificial variables. This result allows us to expand the weighted-sum instances of the DMLC deliverability version of \eqref{eq:propForm} into linear programs for which both extensive theory and efficient methods are available, such as the simplex method and interior point methods (see, e.g., \cite{griva2009} and \cite{wright1997}). The full formulation given in Appendix~\ref{app:propForm} takes this expanded form. A drawback with the expansion is largely increased problem sizes caused by variables and constraints with voxelwise components (usually of the order of $10^5$). In our previous study \cite{engberg2017}, however, we demonstrate how to utilize an interior point method to efficiently solve the deliverability-relaxed version of \eqref{eq:propForm}. The key is to reduce the system of linear equations associated with the interior point method to the order of number of bixels, usually implying a reduction in size by several orders of magnitude. Fortunately, we can obtain an analogous reduction for the DMLC deliverability version, since adding the constraints \eqref{eq:condOfExistence}--\eqref{eq:doseComp} does not change the problem structure. The essential steps to reduce the system of linear equations are summarized below; for linear-algebraic details, we refer to Section~4.A in our previous study. 

The constraint matrix of the linear program expansion is partitioned according to voxelwise dependence; let it be denoted by
\newcommand{\OI}{{\left[\begin{array}{@{\hskip 1pt} 
								    c @{\hskip 1pt}} 
	0 \\ 
    I 
	\end{array}\right]}}
\newcommand{\OIT}{{\left[\begin{array}{@{\hskip 1pt} 
									 c @{\hskip 5pt}
									 c @{\hskip 1pt}} 
	0 & I
	\end{array}\right]}}
\[\arraycolsep=1pt
  \left(\begin{array}{ c c }
        A_{11} & A_{12} \\
        A_{21} & A_{22}
  \end{array}\right),
\]
where $A_{22}$ takes the coefficients for all variables and constraints with voxelwise components. The size of $A_{11}$ mainly depends on the number of bixels, and thus is in general several orders of magnitude smaller than $A_{22}$. The system of linear equations to be solved at each iteration of the interior point method has the form 
\[\arraycolsep=1pt
  \left(\begin{array}{ c @{\hskip -1pt} c @{\hskip 12pt} c @{\hskip 2pt} c }
        -D_1   & 0      & A_{11}^T & A_{21}^T \\[2pt]
        0      & -D_2   & A_{12}^T & A_{22}^T \\[6pt]
        A_{11} & A_{12} & D_3      & 0        \\[2pt]
        A_{21} & A_{22} & 0        & D_4 
  \end{array}\right)
  \Delta x = r_{\text{RHS}},
\]
where $D_1,D_2,D_3$ and $D_4$ are positive definite diagonal matrices containing the so-called complementarity products ($\Delta x$ and $r_{\text{RHS}}$ are generalized unknowns and right-hand-sides). By a simple rearrangement, the voxelwise components are concentrated in the bottom right quadrant:
\[\arraycolsep=1pt
    \left(\begin{array}{ c c @{\hskip 12pt} c c }
        -D_1   & A_{11}^T & 0      & A_{21}^T \\[2pt]
        A_{11} & D_3      & A_{12} & 0        \\[6pt]
        0      & A_{12}^T & -D_2   & A_{22}^T \\ 
        A_{21} & 0        & A_{22} & D_4            
    \end{array}\right)
  \Delta x' = r'_{\text{RHS}}
\]
(the prime indicating the appropriate permutation). Now, for our particular linear program, we have $A_{22} = \begin{bmatrix} 0 & I\end{bmatrix}^T$; the zero and identity parts stem from constraints associated with the maximum/minimum dose and mean-tail-dose objective functions, respectively. This gives a useful structure to the voxelwise quadrant,
\[\arraycolsep=1pt
    \left(\begin{array}{ c c }
        -D_2   & A_{22}^T \\
        A_{22} & D_4
    \end{array}\right)
	= 
	\left(\begin{array}{ c @{\hskip 6pt} c @{\hskip -1pt} c }
        -D_2 & 0       & I       \\[4pt]
        0    & D_{4,1} & 0       \\
        I    & 0       & D_{4,2}
    \end{array}\right), 
    \quad D_4 = 
    \left(\begin{array}{ c @{\hskip -1pt} c }
        D_{4,1} & 0 \\ 
        0 & D_{4,2}
    \end{array}\right)
\]
---but also to its inverse, which is identically structured and can be formed explicitly (see Section~4.A in \cite{engberg2017}). Both solving and multiplying using matrices with this structure are inexpensive operations that merely amount to scaling and adding two vectors. Therefore, the entire system can be solved in two relatively inexpensive steps: by one solve using the Schur complement of the voxelwise quadrant, and by one subsequent operation using the inverse of the voxelwise quadrant. The first step is inexpensive since the Schur complement is considerably smaller than the entire system; its size is identical to that of the top left quadrant, i.e., is of the order of the number of bixels. 

Applying an interior point method to solve the linear program expansion of \eqref{eq:propForm} has an advantage over a simplex method regarding the interpretation and choice of stopping criteria. Once they are feasible, the iterates of the interior point method provide information about the proximity to the optimal value, and the method is usually terminated when the distance is ``sufficiently small''. In the case of a weighted-sum instance of \eqref{eq:propForm}, this stopping criterion is directly interpretable as a dose tolerance (owing to the proposed objective functions being interpretable as dose) and can thus be chosen as such. In our planning study presented in Section~\ref{sec:Results}, the interior point method is terminated when the objective function value is within 1 cGy from the optimal value. While this is to be inexact in a mathematical meaning (the stopping criterion normally relating to machine precision), pushing the method much further is not expected to be meaningful in our application.

\section{Results}\label{sec:Results}
To examine the ability of the proposed objective functions to produce high-quality DMLC plans, we have studied the plan quality indices obtained among DMLC plans Pareto optimal to \eqref{eq:propForm}. DMLC plans generated using the MCO module in RayStation \cite{rs2012} are used as a reference for conventional planning. 

Three patient cases are considered: a prostate, a lung, and a head-and-neck (HN) case. Their respective planning goal amounts to optimizing the three plan quality measures shown in Table~\ref{tab:tradeOffs}, subject to fulfilment of several PTV and OAR dose requirements expressed in terms of dose-at-volume, average, minimum, and maximum dose bounds. Lists of all PTV and OAR requirements are given in Appendix~\ref{app:listOfCriteria}. The restriction to three plan quality indices is to enable the visualization of the obtained values in a three-dimensional coordinate system; in general, any number of indices can be chosen. 
\begin{table}[h]
	\centering
	\begin{minipage}{\textwidth} 
	\renewcommand*\footnoterule{}
	\renewcommand*{\thempfootnote}{\arabic{mpfootnote}}
	\caption{Plan quality indices for the three studied patient cases. \emph{Type} states the type of index, here either the dose-at-volume (d-a-v) at the $v_k$ volume fraction or average dose (avg). For the PTVs, a homogeneity index $\HI_{L\%}^{H\%}$ is constructed as the difference between a low- ($L$) and high-percentage ($H$) dose-at-volume. \emph{Aim} indicates, for clarity, whether the index is subjected to minimization ($\swarrow$) or to maximization ($\nearrow$). Parameters $l_k$ and $u_k$ are hard dose bounds [Gy], while $\hat{l}_k$ and $\hat{u}_k$ are utopian levels [Gy] at which optimizing incentives can be removed.\label{tab:tradeOffs}}
	\begin{tabularx}{\textwidth}{l X X r X r r}
		\toprule
		\vphantom{$\Big($}
		Structure               &                                          & Type  & $v_k$ & Aim        & $l_k,\hat{u}_k$ & $\hat{l}_k,u_k$ \\\midrule
		\vphantom{$\Big($}\it Prostate case                                                                                            \\
		\multirow{2}{*}{PTV 60} & \ldelim\{{2}{0em}[$\HI_{1\%}^{99\%}$\,]  & d-a-v & 99 \% & $\nearrow$ & 57.0, 60.0      & --         \\
		                        &                                          & d-a-v &  1 \% & $\swarrow$ & --              & 60.0, 63.0 \\
		1-cm ring\footnote{1-cm ring structure outside the PTVs} &         & avg   &    -- & $\swarrow$ & --              &  0.0, 60.0 \\
		Rectum                  &                                          & d-a-v & 50 \% & $\swarrow$ & --              &  0.0, 60.0 \\[5pt]
		\vphantom{$\Big($}\it Lung case                                                                                                \\
		\multirow{2}{*}{PTV}    & \ldelim\{{2}{0em}[$\HI_{2\%}^{99\%}$\,]  & d-a-v & 99 \% & $\nearrow$ & 49.5, 55.0      & --         \\
		                        &                                          & d-a-v &  2 \% & $\swarrow$ & --              & 55.0, 60.5 \\
		Lungs                   &                                          & d-a-v & 50 \% & $\swarrow$ & --              &  0.0, 55.0 \\
		Heart                   &                                          & d-a-v & 10 \% & $\swarrow$ & --              &  0.0, 55.0 \\[5pt]
		\vphantom{$\Big($}\it HN case                                                                                                  \\
		\multirow{2}{*}{PTV 56} & \ldelim\{{2}{0em}[$\HI_{20\%}^{95\%}$\,] & d-a-v & 95 \% & $\nearrow$ & 53.2, 56.0      & --         \\
		                        &                                          & d-a-v & 20 \% & $\swarrow$ & --              & 56.0, 69.3 \\
		\multirow{2}{*}{PTV 66} & \ldelim\{{2}{0em}[$\HI_{5\%}^{95\%}$\,]  & d-a-v & 95 \% & $\nearrow$ & 59.4, 66.0      & --         \\
		                        &                                          & d-a-v &  5 \% & $\swarrow$ & --              & 66.0, 69.3 \\
		Parotids                &                                          & avg   &    -- & $\swarrow$ & --              &  0.0, 66.0 \\\bottomrule
	\end{tabularx}
\end{minipage}
\end{table}
Our translation of the planning goals into tricriteria instances of \eqref{eq:propForm} is straightforward: doses-at-volume are replaced with upper and lower mean-tail-dose (average, minimum, and maximum doses are exact, see Appendix~\ref{app:propForm}), their dose bounds left unaltered; and since no optimizing incentive is needed for the PTV and OAR requirements in Appendix~\ref{app:listOfCriteria}, their auxiliary variables are excluded from the objective as described in Section~\ref{sec:Methods}. As to the analogous set-up in RayStation of the tricriteria instances for conventional planning, besides constraints to handle the dose bounds, three penalty-based objective functions are used with the utopian levels $\hat{l}_k$ and $\hat{u}_k$ as reference points. Except for the initial set-up according to the planning goals and constraints given in the tables, no further human-to-TPS interaction is involved in this planning study.

\begin{figure}[h]\centering
	\noindent\makebox[\textwidth]{
	\subfloat[]{\vtop{\vskip0pt\hbox{\centering\includegraphics[scale=.7]{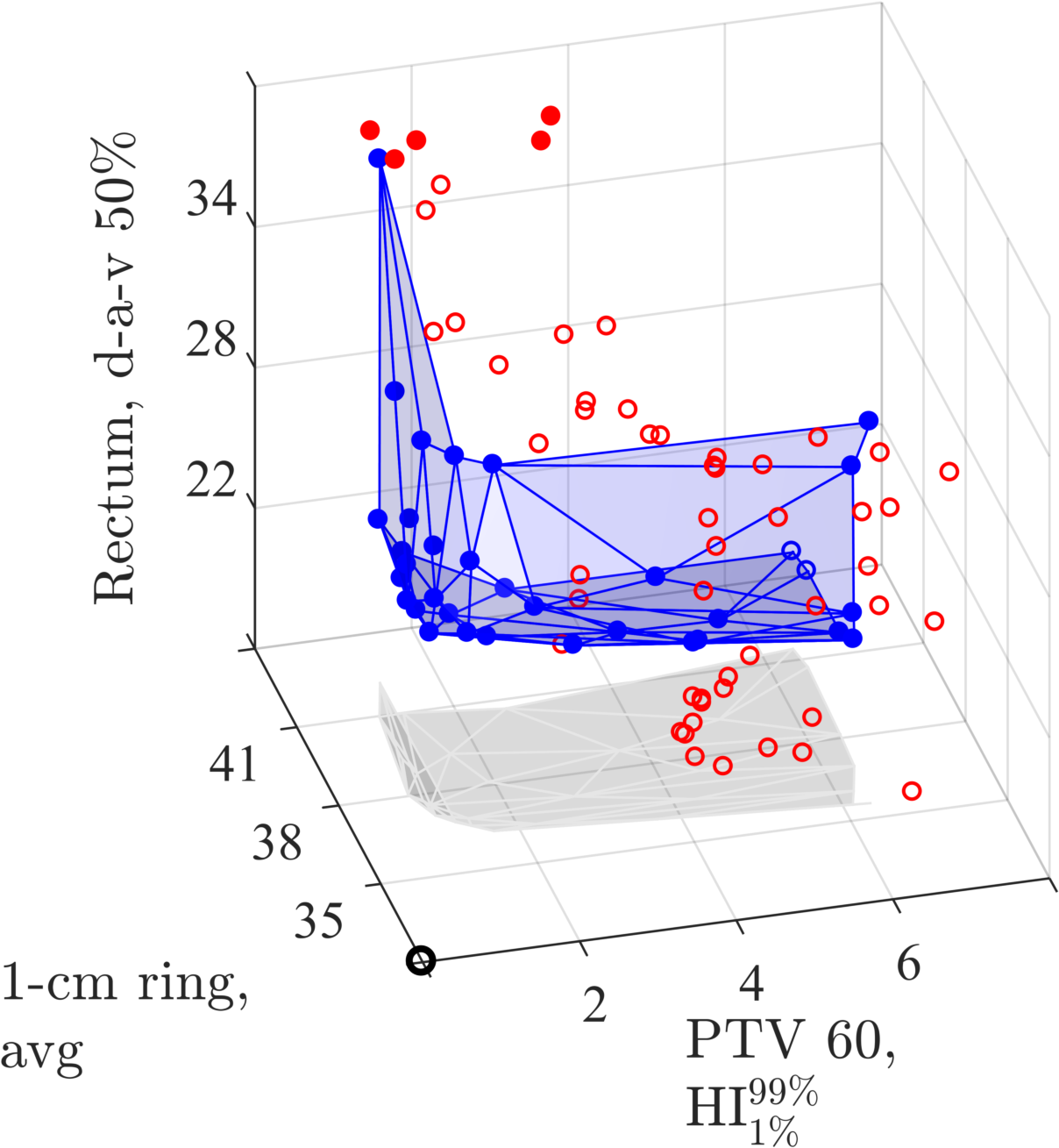}}}}\hspace*{10pt}
	\subfloat[]{\vtop{\vskip0pt\hbox{\centering\includegraphics[scale=.7]{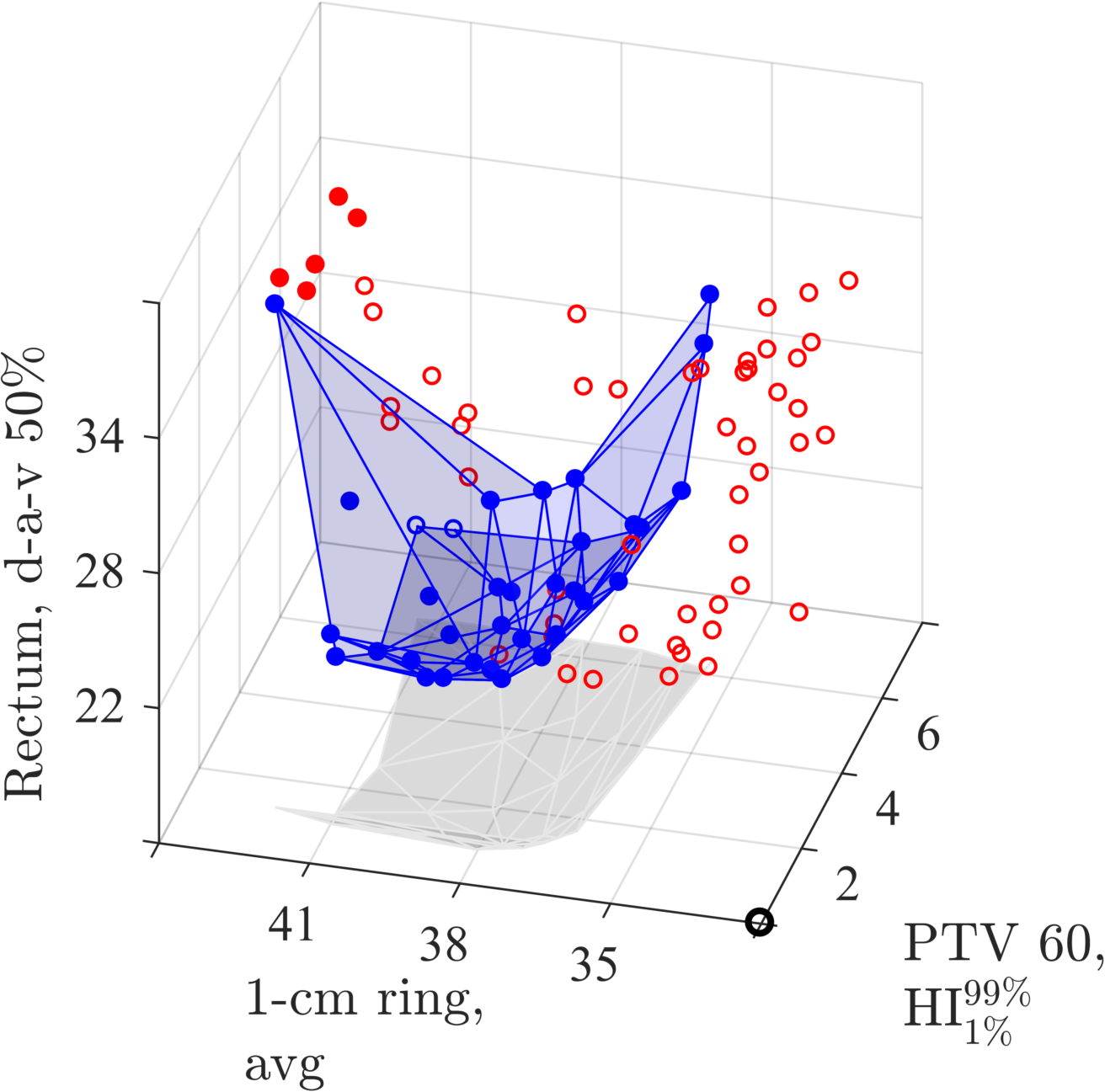}}}}
	}
	\caption{Prostate case. Plan quality indices as listed in Table~\ref{tab:param} obtained among the plans Pareto optimal to \eqref{eq:propForm} (blue dots) with convex hull (blue surface), and obtained among the reference plans Pareto optimal to the conventional formulation (red dots). Unfilled dots correspond to plans that violate PTV constraints by more than 1~\%. The black circle marks the corner of best values, in this case of lowest values. (a) and (b) show different angles of the 3D plot.\label{fig:paretoP}}
\end{figure}

\begin{figure}[h]\centering
	\noindent\makebox[\textwidth]{
	\subfloat[]{\vtop{\vskip0pt\hbox{\centering\includegraphics[scale=.7]{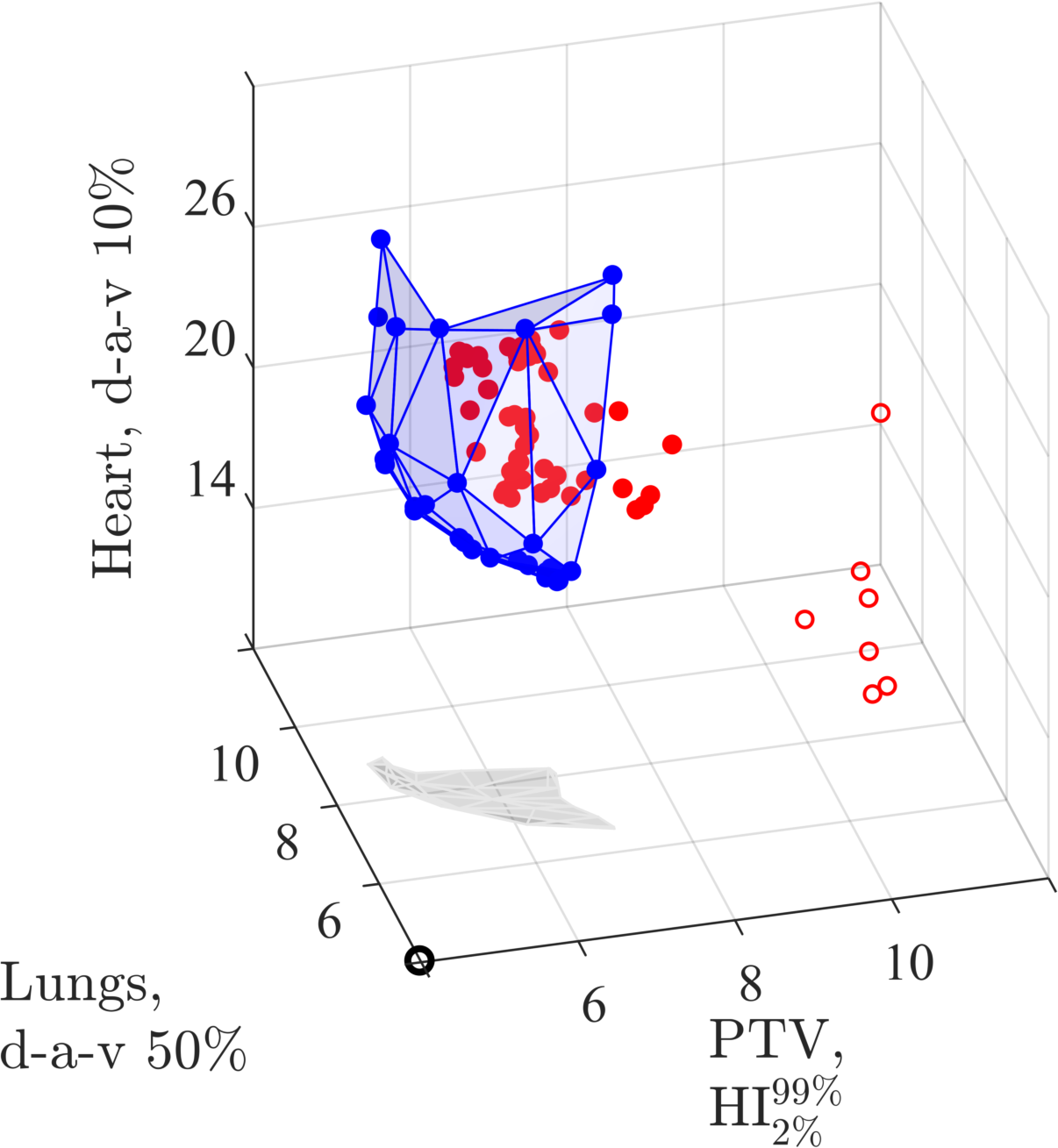}}}}\hspace*{10pt}
	\subfloat[]{\vtop{\vskip0pt\hbox{\centering\includegraphics[scale=.7]{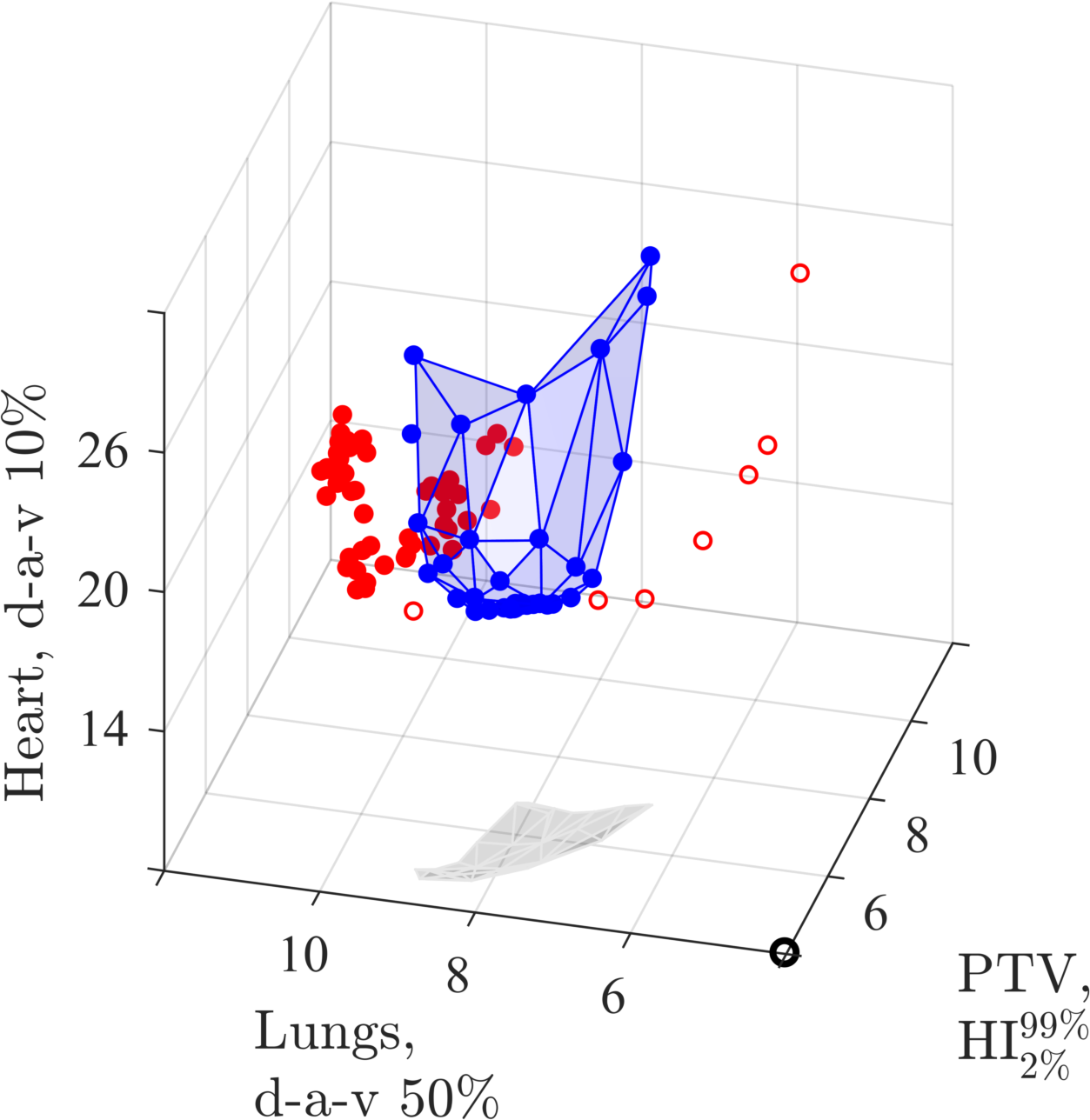}}}}
	}
	\caption{Plan quality indices obtained for the lung case, as listed in Table~\ref{tab:param} and labeled in Figure~\ref{fig:paretoP}. (a) and (b) show different angles of the 3D plot.\label{fig:paretoL}}
\end{figure}

\begin{figure}[h]\centering
	\noindent\makebox[\textwidth]{
	\subfloat[]{\vtop{\vskip0pt\hbox{\centering\includegraphics[scale=.7]{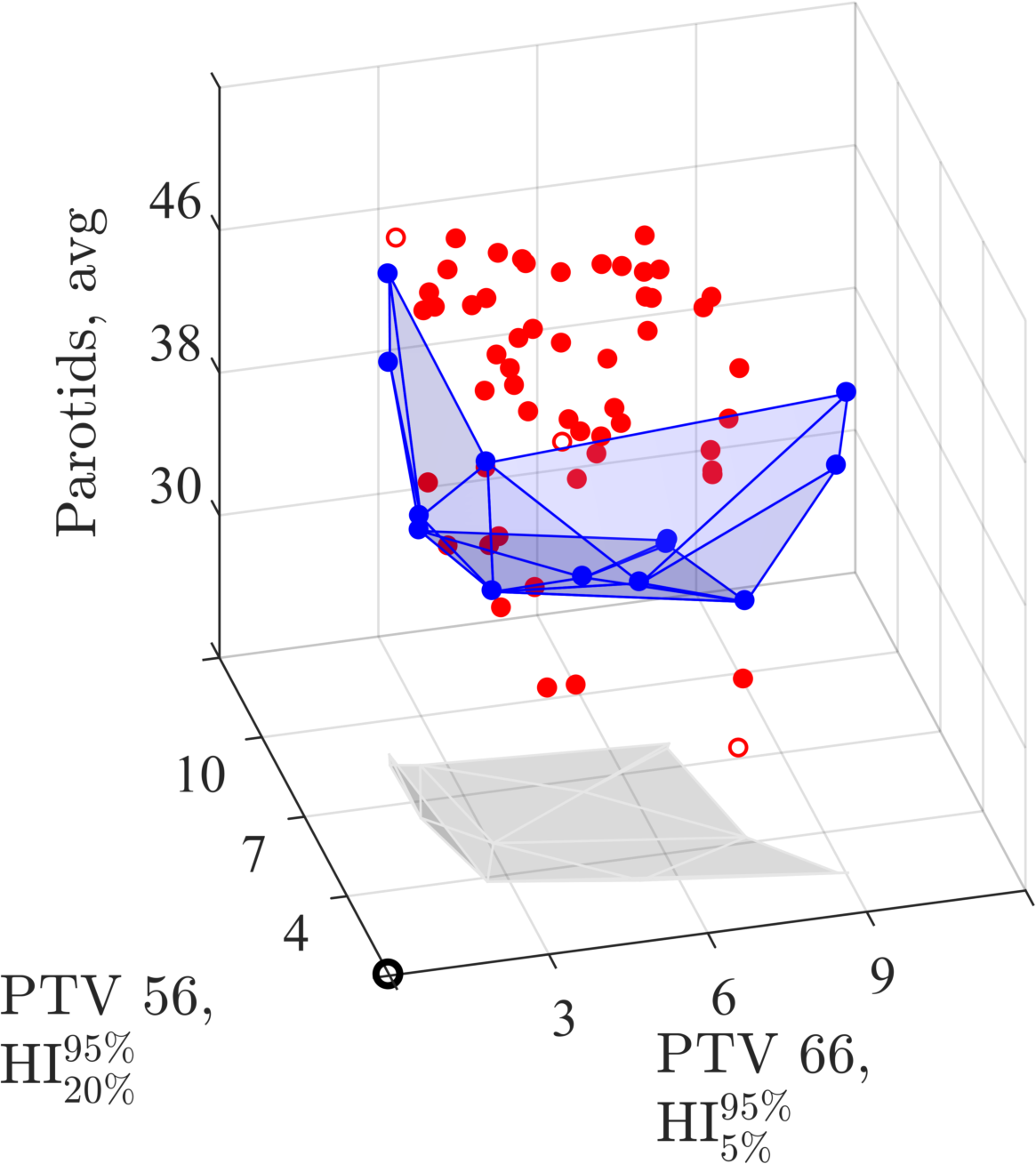}}}}\hspace*{10pt}
	\subfloat[]{\vtop{\vskip0pt\hbox{\centering\includegraphics[scale=.7]{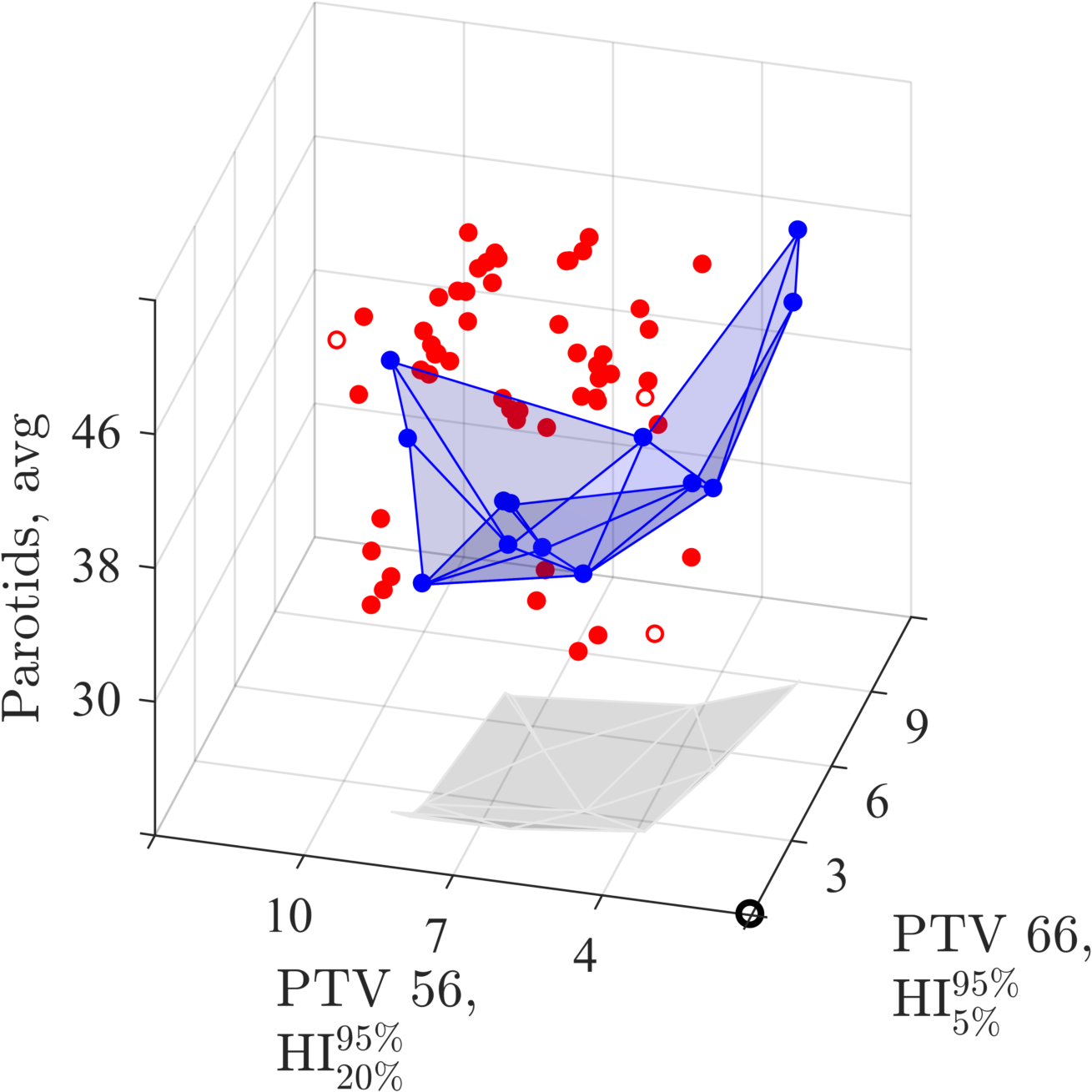}}}}
	}
	\caption{Plan quality indices obtained for the HN case, as listed in Table~\ref{tab:param} and labeled in Figure~\ref{fig:paretoP}. (a) and (b) show different angles of the 3D plot.\label{fig:paretoHN}}
\end{figure}

\begin{figure}[h]\centering
	\noindent\makebox[\textwidth]{
	\subfloat[Prostate case.]{\centering\includegraphics[scale=.7]{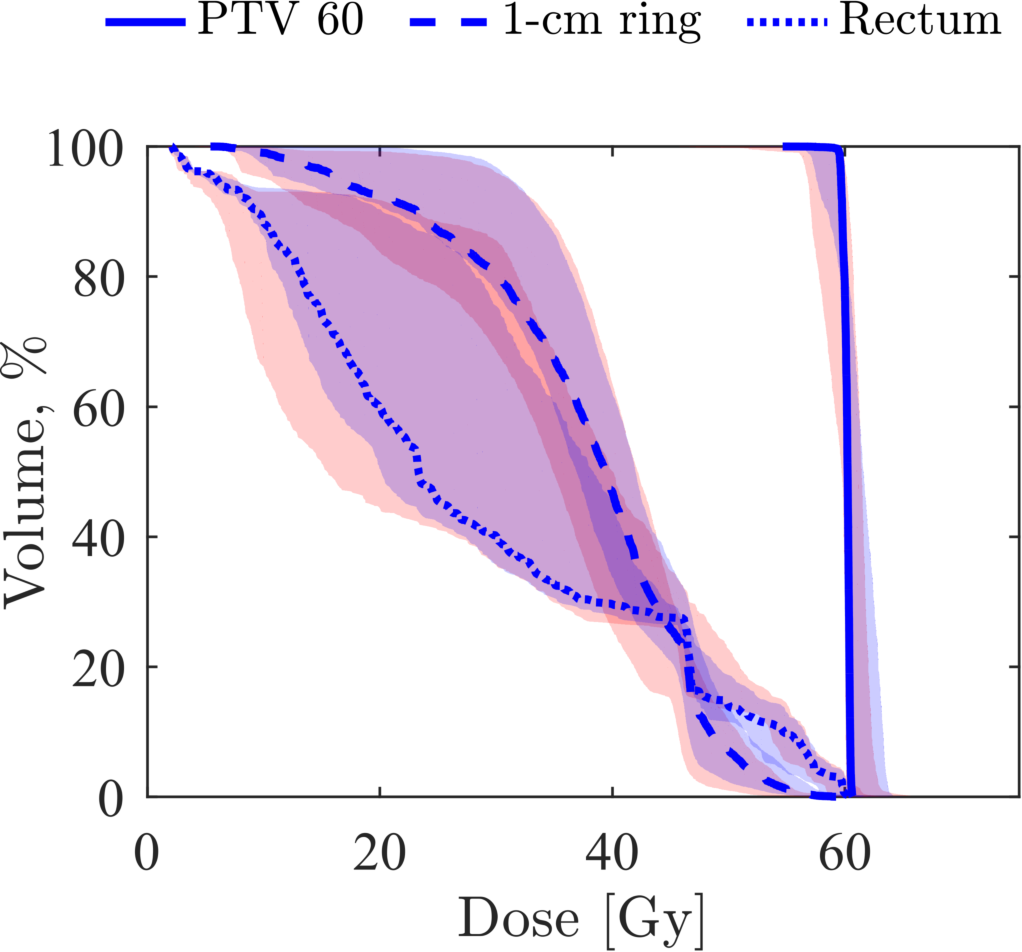}}\hspace*{10pt}
	\subfloat[Lung case.]    {\centering\includegraphics[scale=.7]{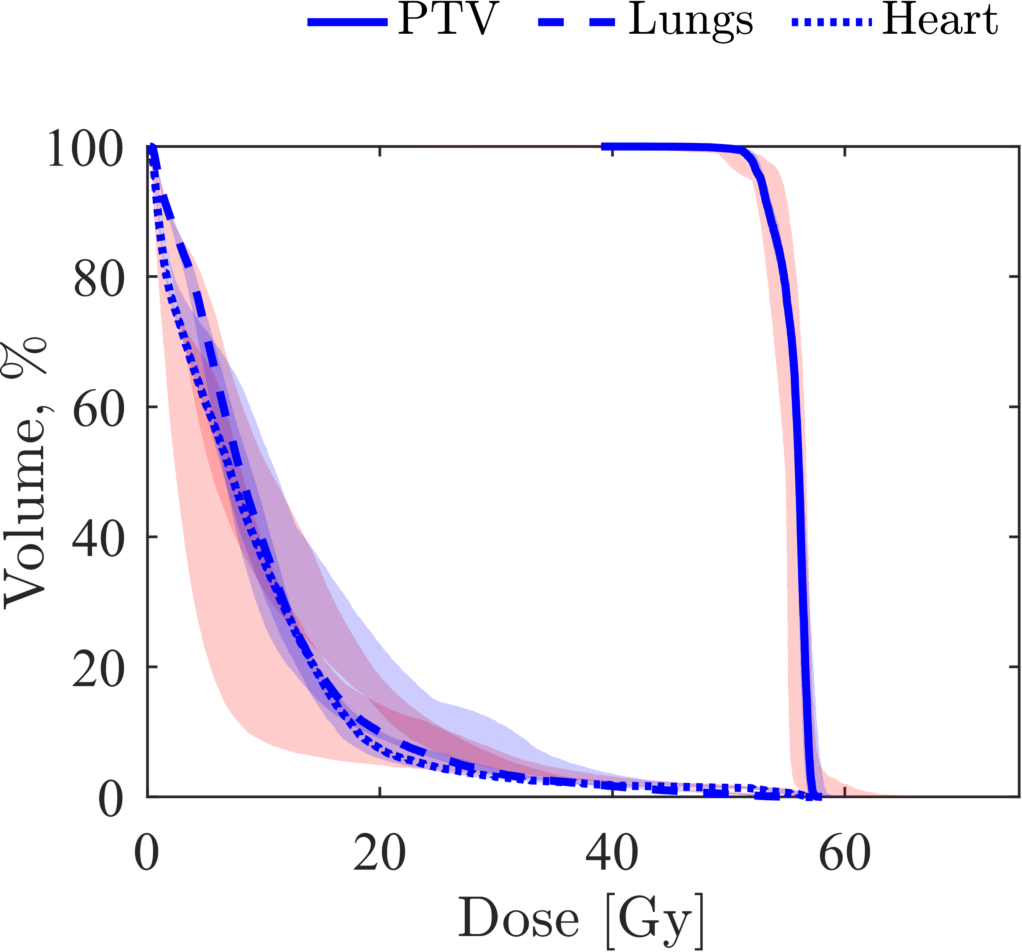}}\hspace*{10pt}
	\subfloat[HN case.]      {\centering\includegraphics[scale=.7]{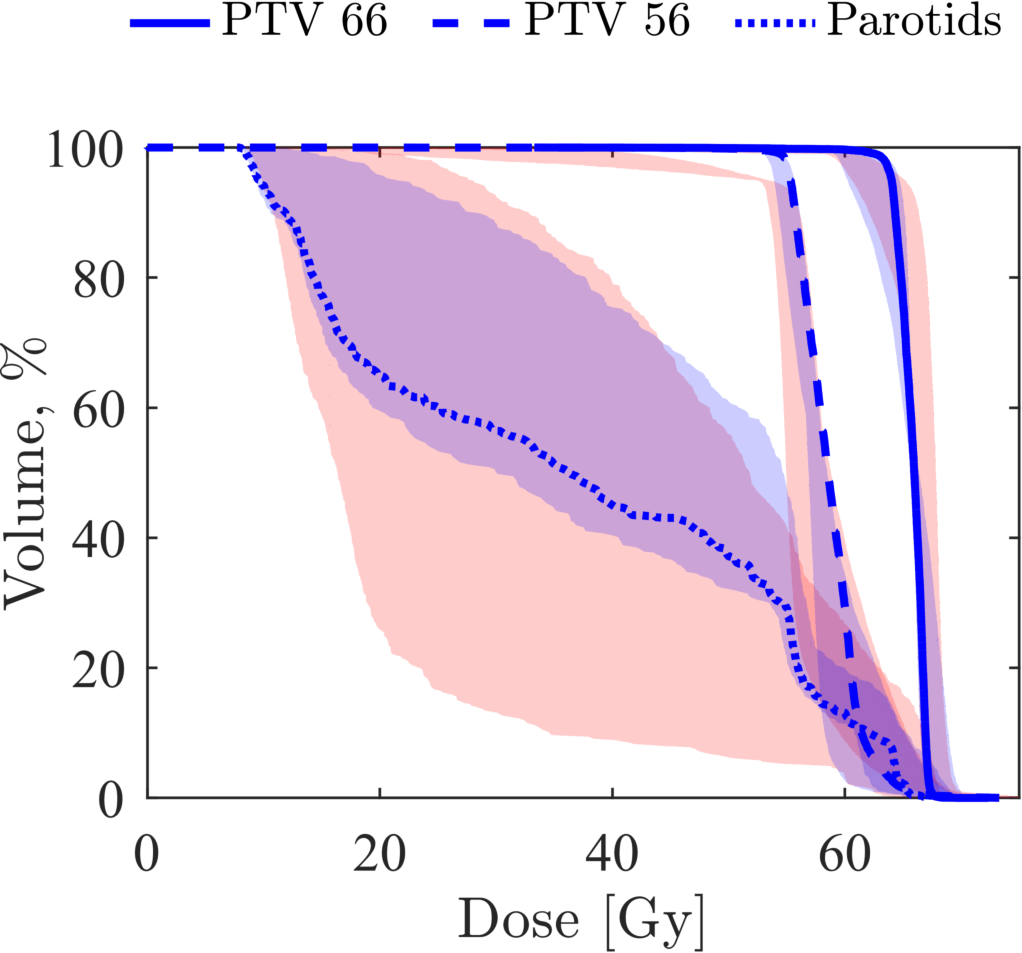}}
	}
	\caption{DVHs for the three patient cases and three respective structures for which plan quality measures are optimized. The blue band covers the area spanned by the plans Pareto optimal to \eqref{eq:propForm}, and the red band covers that spanned by the plans Pareto optimal to the conventional formulation. The highlighted blue histograms correspond to the plan with balanced objective function weights (i.e., with all weights equal to one third).\label{fig:dvhs}}
\end{figure}

\begin{figure}[h]\centering
	\noindent\makebox[\textwidth]{
	\subfloat[Prostate case.]{\vtop{\vskip0pt\hbox{\centering\includegraphics[scale=.7]{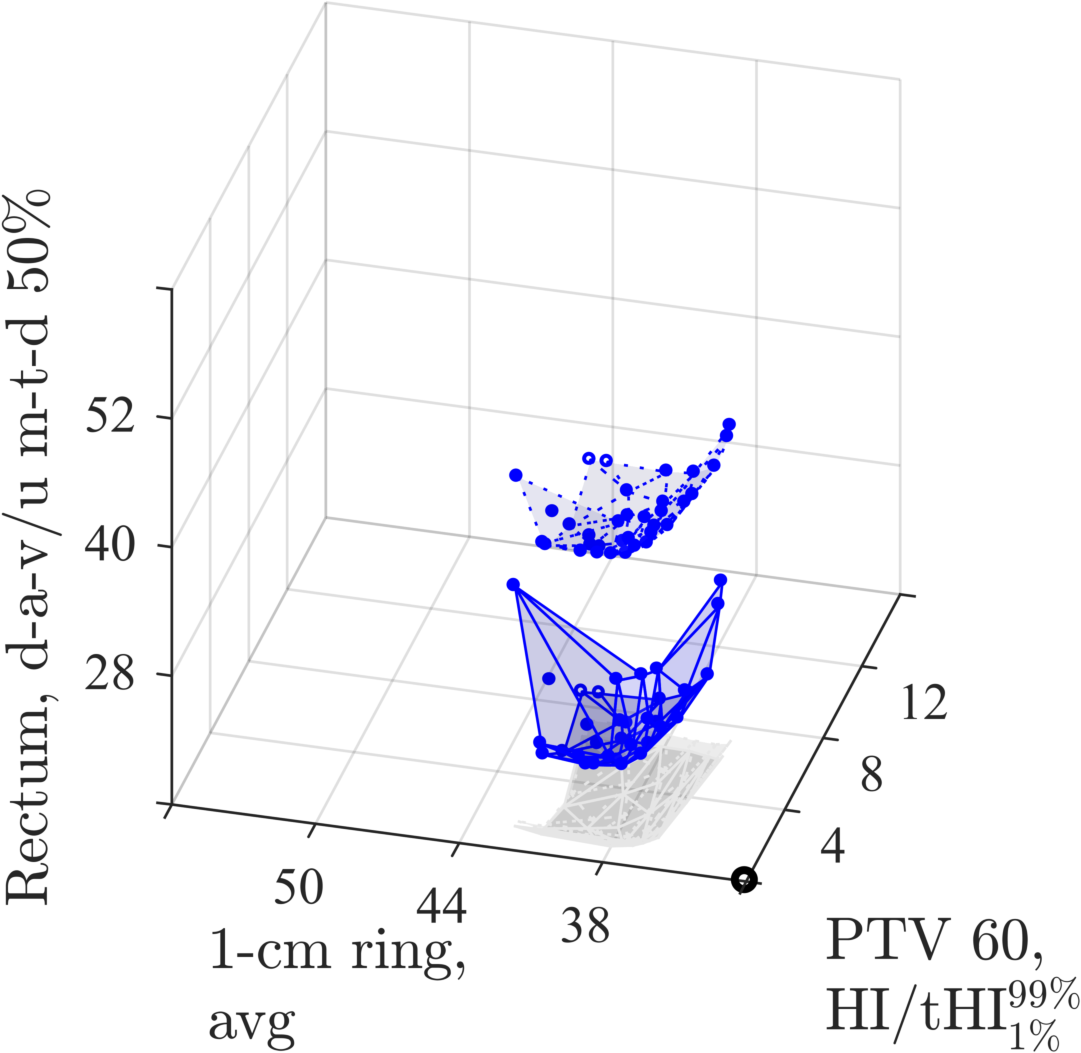}}}}
	\subfloat[Lung case.]{\vtop{\vskip0pt\hbox{\centering\includegraphics[scale=.7]{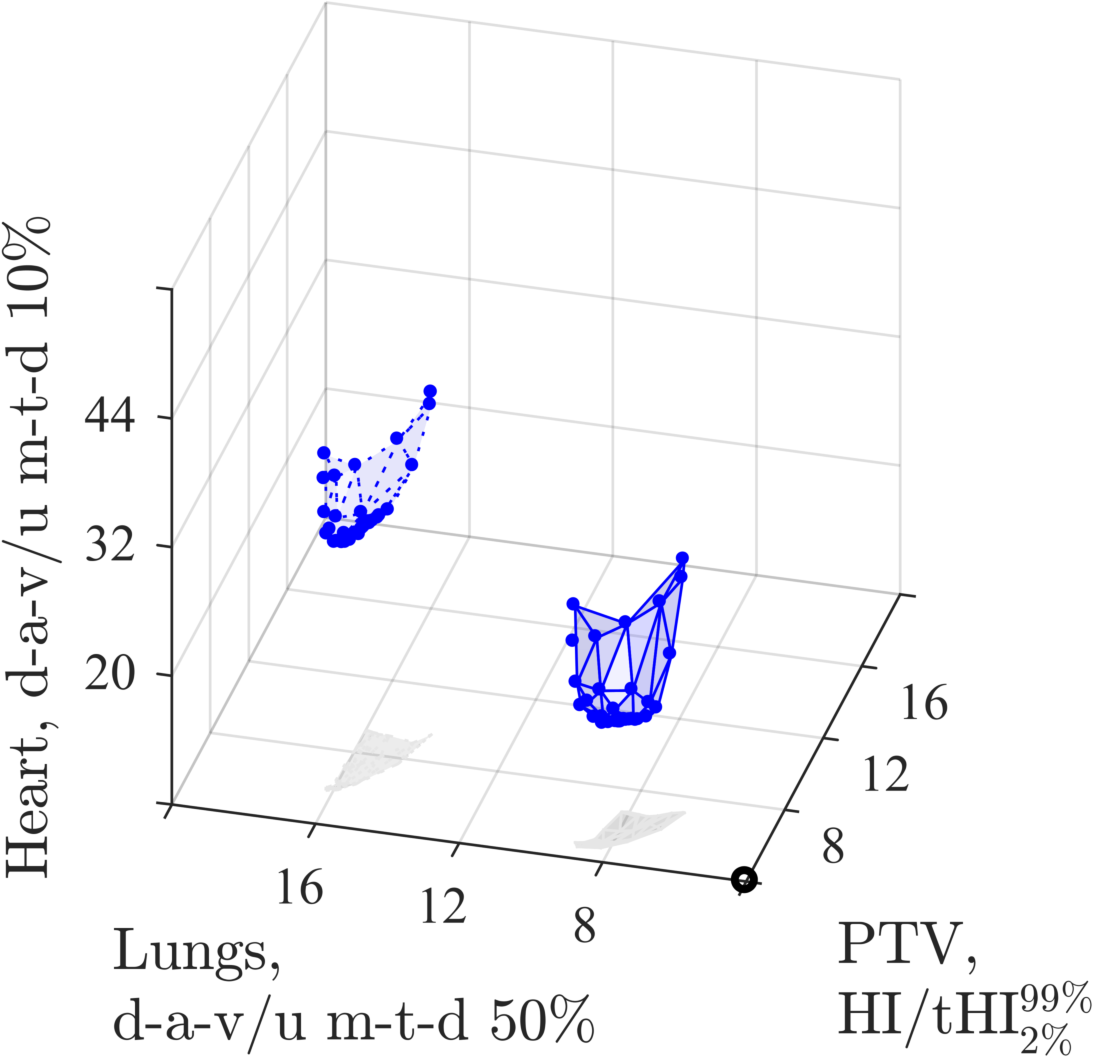}}}}
	\subfloat[HN case.]{\vtop{\vskip0pt\hbox{\centering\includegraphics[scale=.7]{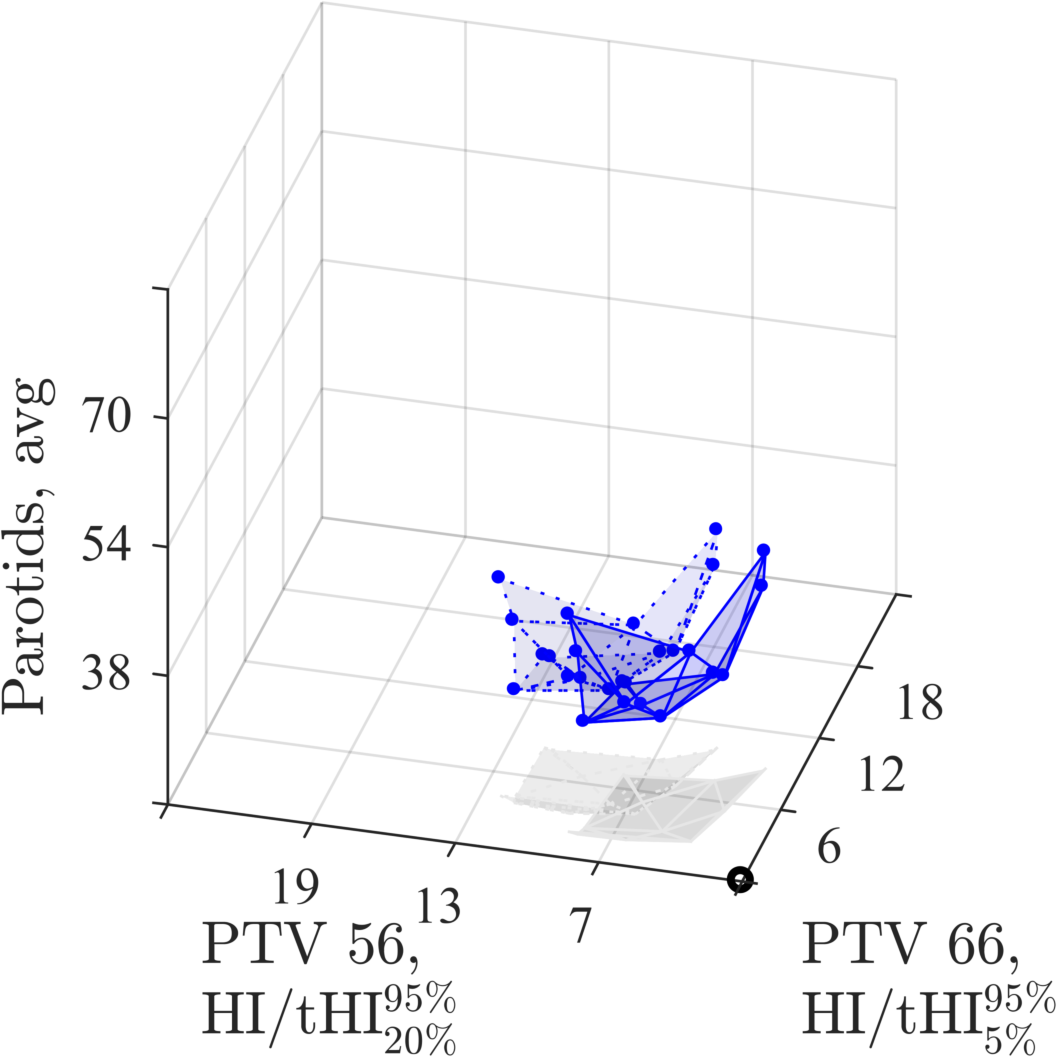}}}}
	}
	\caption{Plan quality indices as labeled in Figures~\ref{fig:paretoP}--\ref{fig:paretoHN} (solid convex hull) in juxtaposition with the corresponding Pareto optimal objective function values of \eqref{eq:propForm} (dotted convex hull), i.e., upper mean-tail-doses (u m-t-d), average doses (avg), and the mean-tail-dose analogue of the homogeneity index ($\tHI^{H\%}_{L\%}$) (the difference between a low-percentage ($L$) upper, and a high-percentage ($H$) lower mean-tail-dose).\label{fig:paretoCVaR}}
\end{figure}

A total of 34 DMLC plans Pareto optimal to \eqref{eq:propForm} are generated for the prostate and lung case, and a total of 13 for the computationally more expensive HN case. The underlying objective function weights are normalized and equidistantly distributed on the triangle with corners in (1,0,0), (0,1,0), and (0,0,1). As to the reference DMLC plans Pareto optimal to the conventional formulation, a total of 55 are generated for each case using the MCO module in RayStation. Independently of how generated, all plans undergo final dose computation of clinical accuracy and are then evaluated based on the three targeted plan quality indices of Table~\ref{tab:tradeOffs}. The outcome of this evaluation is shown in Figures~\ref{fig:paretoP}--\ref{fig:paretoHN}, and DVHs are shown in Figure~\ref{fig:dvhs} to provide a broader perspective on the plans with regards to dose distribution. Figure~\ref{fig:paretoCVaR} shows the plan quality indices in juxtaposition with the Pareto optimal objective function values of \eqref{eq:propForm}.

Figures~{\ref{fig:paretoP}}-{\ref{fig:paretoHN}} play the central role in our analysis of the results, which are similar for all three patient cases. In these figures, the generated plans are represented by their plan quality indices---i.e., \emph{not} by their respective Pareto optimal objective function values. As expected from the weak connection between the penalty-based objective functions and the plan quality measures, the plans Pareto optimal to the conventional formulation (red dots) appear arbitrarily scattered with this representation. In contrast, the plans Pareto optimal to \eqref{eq:propForm} (blue dots) conform to a non-dominated surface in a similar way as Pareto optimal values do to a Pareto surface (we refer to a point as dominating over another point if it dominates in all coordinates). This conformity is indicated by the fact that almost all of these plans lie on the convex hull (blue surface). In fact, as can be observed in Figure~\ref{fig:paretoCVaR}, the convex hull of the plan quality indices appears only rigidly shifted relative to the convex hull of the corresponding Pareto optimal values, rather than deformed in shape. These indications support the assumed stronger correlation between the proposed objective functions and the plan quality measures.

As to the values obtained in Figures~{\ref{fig:paretoP}}-{\ref{fig:paretoHN}}, none of the two sets of points appears to in a clear way dominate the other. Dominance is particularly unclear for the prostate case shown in Figure~\ref{fig:paretoP}. However, the situation changes when violation of constraints is taken into account; violation of PTV constraints by more than 1~\% is marked in the figures by unfilled dots. A dominance in favor of the proposed objective functions is revealed when violating plans are ignored. For example, for the lung case shown in Figure~\ref{fig:paretoL}, the violating plans are exactly those out of the conventionally generated plans that are not dominated. Although minor violation is somewhat expected for all plans due to the final accurate dose computation, the violation obtained among plans generated using the conventional objective functions is more systematic and of a larger extent. This inherent tendency of using the conventional approach, giving solutions that sometimes fail to satisfy constraints, is addressed and explained by Fredriksson \cite{fredriksson2012}, who suggests using smoother penalty functions as a possible remedy.

Implications of the results are further discussed in Section~\ref{sec:Discussion}.

\section{Discussion}\label{sec:Discussion}
The goal of this study was to examine the streamlining potential on the planning process of the objective functions suggested in our previous study \cite{engberg2017} and formulated in \eqref{eq:propForm} when stress-tested with DMLC deliverability constraints and accurate dose computations. The proposed objective functions are designed to offer a stronger connection to measures of plan quality than the conventional penalty-based objective functions, and to give tractable (convex) optimization problems. The proposed functions are argued to have a streamlining potential if they, in addition, are able to produce plans of at least equivalent quality as the conventional functions. The outcome of our planning study indeed indicates such an ability. 

Surprisingly, the results presented in Figures~\ref{fig:paretoP}--\ref{fig:paretoHN} and commented in Section~\ref{sec:Results} were obtained despite several handicaps for the proposed formulation. One such is the already addressed unfair advantage for the conventional formulation of systematically violated constraints. Violated constraints theoretically give more room to optimize the objective function values (here, to minimize the penalties), yet, the plan quality indices obtained are comparable. Another handicap is the fact that the dose bounds aimed for PTV and OAR doses-at-volume are imposed unaltered on the mean-tail-doses when optimizing using the proposed formulation. Since upper and lower mean-tail-doses over- and underestimate dose-at-volume, respectively, the unaltered bounds overly restrict (i.e., narrow) the feasible region to the optimization problem. A smaller feasible region gives less room to optimize the objective function values, yet again, none of the two sets of plan quality indices dominates the other. It seems indeed that the characterizing properties of the proposed formulation---consistency with respect to plan quality measures, and convexity---are able to compensate for these handicaps.

As mentioned in Section~\ref{sec:Results}, no human-to-TPS interaction was involved in our planning study except for initial set-up. The results are the outcome of optimizing an MCO set-up based on a straightforward and systematic translation of the clinical goals and requirements into objective functions and constraints. Re-optimizing each of the conventionally generated plans with revised parameters (i.e., revised reference points) could improve their quality and feasibility. The extent of this potential improvement has not been investigated, but what \emph{can} be concluded is an existing need to re-optimize all three studied patient cases in order to meet the outcome from using \eqref{eq:propForm}.

The positive effects of explicitness and convexity of the proposed formulation \eqref{eq:propForm} do not come for free. In abandoning the penalty function paradigm in favor of these aspects, compromises have been made with regards to computational complexity. The problem size of \eqref{eq:propForm} in its expanded form is proportional to the number of voxels, which can be contrasted to the problem size of the conventional formulation proportional to the usually substantially smaller number of bixels. Although fast optimization possibly becomes less important with more accurate planning tools, running times must be reasonable to justify practical use. We have demonstrated in Section~\ref{sec:Methods:Solvability} how to exploit the structure of \eqref{eq:propForm} to significantly reduce the computational cost associated with solving this formulation using a tailored interior point method. Combined with the techniques developed by Colombo and Gondzio \cite{colombo2008} to further reduce this cost, running times of the order of 5, 10, and 60 minutes per weighted-sum instance of the prostate, lung, and HN case were required to meet the desired optimization tolerance with our MATLAB (MathWorks, Natick, Massachusetts) implementation on an Intel Core i7 2.80 GHz computer with four cores. We anticipate a decrease if GPUs are used. 

Treatment planning is a complex task. To what extent it is complicated, however, depends on the planning tools available to support the planner in completing that task. This distinction between the \emph{complexity} of a process, and \emph{complicated tools} to handle the process is in line with the theories and experiences presented by Andersson et al.~\cite{andersson2014}. While the latter is a consequence of shortcomings in tool design and can (should) be avoided, complexity is something inherent that needs to be accepted. In our research, we embrace the complexity of treatment planning while aiming at making the planning process less complicated, specifically, by offering more accurate tools, i.e., by reformulating the underlying objective functions to more comprehensive alternatives. Another example of reducing complication while accepting the complexity of treatment planning, is the successful \cite{craft2012} clinical introduction of a posteriori MCO techniques made possible through rigorous and extensive research (see, e.g., \cite{bokrantz2013,craft2006,monz2008}). With these techniques at hand, the otherwise tedious manual task of finding weighting factors that give a desired objective function trade-off can be accomplished in a formalized and structured manner without hiding degrees of freedom; the weighting factors indeed introduce a flexibility to the planning process that should be available for the planner to explore. We believe that combining MCO techniques with the proposed objective functions can provide the skilled planner with the necessary tools to control the planning process in an uncomplicated fashion. Future studies should be aiming at further developing practical aspects of using the proposed formulation, in particular with regards to compliance with other treatment delivery techniques and with other types of plan quality indices.

\section{Conclusion}
We have investigated the ability of our previously suggested formulation of objective functions for treatment plan MCO to produce DMLC plans of plan quality comparable to that of DMLC plans generated using the conventional penalty-based objective functions. In a numerical planning study involving three patient case, DMLC plans Pareto optimal to the proposed formulation were generated and compared to the DMLC plans resultant from the penalty-based MCO module in RayStation. Comparable plan quality---or, better when accounting for the systematic violation of constraints among conventionally generated plans---was observed when evaluating the plans based on three pre-defined plan quality indices. Supported by these results, the proposed objective functions are argued to have a potential to streamline the planning process, since they also are designed to overcome methodological shortcomings in the conventional objective functions that cause a need for trial and error and time-consuming re-optimizations in the current planning process. These conclusions position the proposed formulation as an alternative to existing methods for automated planning. In particular, we believe that the increased accuracy of the planning tools imposed by the proposed objective functions is essential in making the planning process less complicated.

\bibliography{dmlc}
\bibliographystyle{siam} 

\raggedbottom
\pagebreak

\appendix

\section{Full proposed formulation}\label{app:propForm}
In the full version of \eqref{eq:propForm}, we include objective functions to minimize ($\swarrow$) dose-at-volume (d-a-v) (via mean-tail-dose approximation), maximum dose (max), and average dose (avg) (given by \eqref{eq:fullPropFormUMTD-a}--\eqref{eq:fullPropFormUMTD-b}, \eqref{eq:fullPropFormMAX}, and \eqref{eq:fullPropFormUAVG}), as well as to maximize ($\nearrow$) dose-at-volume (via lower mean-tail-dose approximation), minimum dose (min), and average dose (\eqref{eq:fullPropFormLMTD-a}--\eqref{eq:fullPropFormLMTD-b}, \eqref{eq:fullPropFormMIN}, and \eqref{eq:fullPropFormLAVG}). The objective functions, collected according to type and aim in index sets $K_{\text{type}}^{\text{aim}}$, are assumed sorted in such a way that 
\[
	K_{\text{d-a-v}}^{\swarrow} \cup K_{\max}^{\swarrow} \cup K_{\text{avg}}^{\swarrow} = \left\{1,\ldots,q\right\}
\]
and
\[
	K_{\text{d-a-v}}^{\nearrow} \cup K_{\min}^{\nearrow} \cup K_{\text{avg}}^{\nearrow} = \left\{q\!+\!1,\ldots,K\right\}. 
\]
The full proposed formulation is then given by
\begin{subequations}\label{eq:fullPropForm}
\begin{alignat}{3}
& \minimize{\alpha_k,\xi_k,\eta^k,d} 
           && \quad \mkern-10mu \mathmakebox[0pt][l]{\big[\,\xi_1, \cdots, \xi_q, -\xi_{q+1}, \cdots, -\xi_K\,\big]^T} \nonumber\\
& \subject && \quad \alpha_k + \frac{1}{v_k}\sum_{i \in V_{s_k}} \Delta^{s_k}_i \eta^k_i \leq \xi_k,
                                                           && k \in K_{\text{d-a-v}}^{\swarrow}, \label{eq:fullPropFormUMTD-a}\\
&          && \quad \eta^k_i \geq d_i - \alpha_k, \enskip \eta^k_i \geq 0, \enskip i \in V_{s_k}, \enskip  
                                                           && k \in K_{\text{d-a-v}}^{\swarrow}, \label{eq:fullPropFormUMTD-b}\\[6pt]
&          && \quad d_i \leq \xi_k, \enskip i \in V_{s_k}, && k \in K_{\max}^{\swarrow},         \label{eq:fullPropFormMAX}   \\[6pt]
&          && \quad \sum_{i \in V_{s_k}} d_i \leq \xi_k,   && k \in K_{\text{avg}}^{\swarrow},   \label{eq:fullPropFormUAVG}  \\[2pt]
&          && \quad \hat{l}_k \leq \xi_k \leq u_k,         && k = 1,\ldots,q,                                        \nonumber\\[8pt]
&          && \quad \alpha_k - \frac{1}{1-v_k}\sum_{i \in V_{s_k}} \Delta^{s_k}_i \eta^k_i \geq \xi_k, 
                                                           && k \in K_{\text{d-a-v}}^{\nearrow}, \label{eq:fullPropFormLMTD-a}\\
&          && \quad \eta^k_i \geq \alpha_k - d_i, \enskip \eta^k_i \geq 0, \enskip i \in V_{s_k}, 
                                                           && k \in K_{\text{d-a-v}}^{\nearrow}, \label{eq:fullPropFormLMTD-b}\\[6pt]
&          && \quad d_i \geq \xi_k, \enskip i \in V_{s_k}, && k \in K_{\min}^{\nearrow},         \label{eq:fullPropFormMIN}   \\[6pt]
&          && \quad \sum_{i \in V_{s_k}} d_i \geq \xi_k,   && k \in K_{\text{avg}}^{\nearrow},   \label{eq:fullPropFormLAVG}  \\[2pt]
&          && \quad \hat{u}_k \geq \xi_k \geq l_k,         && k = q\!+\!1,\ldots,K,                                  \nonumber\\[8pt]
&          && \quad \mathmakebox[0pt][l]{d \text{ deliverable by linear constraints \eqref{eq:condOfExistence}--\eqref{eq:doseComp},}}\nonumber
\end{alignat}
\end{subequations}
where $\alpha_k$ and $\eta^k$ are artificial variables used to linearly handle the mean-tail-dose functions. The dose distribution $d$, auxiliary variables $\xi_k$, volume fraction $v_k$, and bounds $l_k,\hat{l}_k,u_k,\hat{u}_k$ are as defined in Section~\ref{sec:Methods}. Moreover, $V_s$ collects the voxels of ROI $s$, and $\Delta^s_i$ is the relative volume of $s$ located in voxel $i$ such that $\sum_{i \in V_s} \Delta^s_i = 1$. Notice that the constraints \eqref{eq:fullPropFormUMTD-b}, \eqref{eq:fullPropFormMAX}, \eqref{eq:fullPropFormLMTD-b}, and \eqref{eq:fullPropFormMIN} and the variable $\eta^k$ have voxelwise components, as addressed in Section~\ref{sec:Methods:Solvability}. 

\section{List of PTV and OAR requirements}\label{app:listOfCriteria}
Table~\ref{tab:param} shows the PTV and OAR requirements for the prostate, lung and HN case introduced in Section~\ref{sec:Results}. In addition to these constraints, dose fall-off is required in terms of average and maximum dose bounds for several ring structures at different distances from the PTVs. 
\begin{table}[!h]
	\centering
	\caption{PTV and OAR dose requirements. \emph{Type} indicates the type of criteria, here dose-at-volume (d-a-v) at the $v_k$ volume fraction, or minimum (min), maximum (max), or average (avg) dose. Parameters $l_k$ and $u_k$ [Gy] give the dose bounds.\label{tab:param}}
	\noindent\makebox[\textwidth]{
		\begin{minipage}{.55\textwidth}
			\subfloat[Prostate case criteria.\label{tab:paramP}]{
				\begin{tabularx}{\textwidth}{l X r r r}
					\toprule
					\vphantom{$\Big($}
					Structure   & Type  & $v_k$ & $l_k$ & $u_k$ \\\midrule
					Bladder     & d-a-v &  5 \% &    -- &  60.0 \\
					            & d-a-v & 25 \% &    -- &  48.6 \\
					            & d-a-v & 50 \% &    -- &  40.8 \\
					Hips        & d-a-v & 50 \% &    -- &  40.8 \\
					Penile bulb & d-a-v & 10 \% &    -- &  48.6 \\
					            & d-a-v & 50 \% &    -- &  40.8 \\ 
					Rectum      & d-a-v &  3 \% &    -- &  60.0 \\
					            & d-a-v & 15 \% &    -- &  57.0 \\
					            & d-a-v & 30 \% &    -- &  52.8 \\ 
					PTV 60      & d-a-v & 50 \% &  60.0 &    -- \\
					PTV 48-60   & d-a-v & 50 \% &  48.0 &    -- \\
					PTV 48      & d-a-v & 99 \% &  45.6 &    -- \\
					            & max   &    -- &    -- &  66.0 \\\bottomrule
				\end{tabularx}
			}\\
			\subfloat[Lung case criteria.\label{tab:paramL}]{
				\begin{tabularx}{\textwidth}{l X r r r}
					\toprule
					\vphantom{$\Big($}
					Structure    & Type  & $v_k$ & $l_k$ & $u_k$ \\\midrule
					Brac. plexus & d-a-v &  1 cc &         -- &       53.0 \\
					Heart        & d-a-v & 33 \% &         -- &       40.0 \\
					Lungs        & avg   &    -- &         -- &       15.0 \\
					             & d-a-v & 20 \% &         -- &       25.0 \\
					             & d-a-v & 30 \% &         -- &       20.0 \\
					Oesophagus   & d-a-v & 10 cc &         -- &       55.0 \\ 
					Cord, 5 mm   & d-a-v &  1 cc &         -- &       44.0 \\
					Cord         & d-a-v &  1 cc &         -- &       42.0 \\
					PTV          & d-a-v & 95 \% &       52.3 &         -- \\
					             & d-a-v & 50 \% &       55.0 &         -- \\
					             & d-a-v &  5 \% &         -- &       57.8 \\\bottomrule
				\end{tabularx}
			}
		\end{minipage}
		\quad
		\begin{minipage}{.62\textwidth}\vspace*{-50pt}
			\renewcommand*\footnoterule{}
			\renewcommand*{\thempfootnote}{\arabic{mpfootnote}}
			\subfloat[HN case criteria.\label{tab:paramHN}]{
				\begin{tabularx}{\textwidth}{l X r r r}
					\toprule
					\vphantom{$\Big($}
					Structure       & Type  & $v_k$ & $l_k$ & $u_k$ \\\midrule
					Pacemaker       & max   &    -- &    -- &   3.5 \\
					Cord            & max   &    -- &    -- &  46.0 \\
					                & d-a-v &  1 cc &    -- &  44.0 \\
					Cord, 5 mm      & max   &    -- &    -- &  48.0 \\
					                & d-a-v &  1 cc &    -- &  46.0 \\
					Brainstem       & max   &    -- &    -- &  54.0 \\ 
					                & d-a-v &  1 cc &    -- &  52.0 \\
					Brainstem, 5 mm & d-a-v &  1 cc &    -- &  54.0 \\
					Brain           & max   &    -- &    -- &  60.0 \\
					                & d-a-v &  1 cc &    -- &  58.0 \\
					Chiasm          & max   &    -- &    -- &  48.0 \\
					Chiasm, 3 mm    & max   &    -- &    -- &  50.0 \\
					Opt. nerves     & max   &    -- &    -- &  48.0 \\
					Opt. nerves, 3 mm & max &    -- &    -- &  50.0 \\
					Orbits          & max   &    -- &    -- &  40.0 \\
					Lenses          & max   &    -- &    -- &   6.0 \\
					Mandible        & d-a-v &  1 cc &    -- &  69.3 \\
					PTV 56          & d-a-v & 99\footnote{relaxed to 95 \% for conventionally generated plans} \% &  53.2 &    -- \\
					                & avg   &    -- &    -- &  58.8 \\
					                & d-a-v &  1 \% &    -- &  69.3 \\
					PTV 66          & d-a-v & 99 \% &  59.4 &    -- \\
					                & avg   &    -- &  66.0 &    -- \\
					                & d-a-v &  1 \% &    -- &  69.3 \\\bottomrule
				\end{tabularx}
			}
		\end{minipage}
	}
\end{table}

\end{document}